\documentclass[12pt, a4paper]{amsart}

\usepackage[all]{xy}
\usepackage{amsmath}
\usepackage{amssymb}
\usepackage{amsfonts}
\usepackage{graphicx}
\usepackage[mathscr]{eucal}
\usepackage{enumerate}
\usepackage[latin1]{inputenc}

\theoremstyle{definition}
\newtheorem{definition}{Definition}[section]
\theoremstyle{remark}
\newtheorem*{rem}{\textbf{Remark}}

\newtheoremstyle{estilo}{\topsep}{\topsep}{\slshape}{}{\bfseries}{.}{ }{}
\theoremstyle{estilo}
\newtheorem{theorem}[definition]{Theorem}
\newtheorem{corollary}[definition]{Corollary}
\newtheorem{lemma}[definition]{Lemma}

\def\ot{\otimes}
\def\op{\oplus}

\def\cd{\cdot}

\renewenvironment{proof}{{\noindent\sc Proof\;}}{\qed\\}

\DeclareMathOperator{\Id}{Id} \DeclareMathOperator{\End}{End}

\def\O{\Omega}

\newcommand{\dtext}[1]{\emph{\textbf{#1}}}
\def\vphi{\varphi}
\def\eps{\varepsilon}
\newcommand{\abs}[1]{\left\lvert#1\right\rvert}

\def\n{\mathbb N}

\def\l{\lambda}
\def\lto{\longrightarrow}
\def\lmto{\longmapsto}

\def\cam{\mathfrak{X}}

\title{Connections over twisted tensor products of algebras}
\author{Javier L\'{o}pez Pe\~na}
\address{Department of Algebra,
University of Granada\\
Avda. Fuentenueva s/n, E-18071, Granada, Spain}
\email{jlopez@ugr.es}
\thanks{This research has been partially supported
by the projects MTM2004-08125 and FQM-266 (Junta de Andaluc\'{i}a
Research Group) and Spanish MEC-FPU grant AP2003-4340. The author would like to
thank J.A. Ca\~nizo for
enlightening discussions.}
\date{}

\begin{document}

\allowdisplaybreaks
\begin{abstract}
    Motivated by some results in classical differential geometry, we give a
constructive procedure for building up a connection over a (twisted) tensor
product of two algebras, starting from connections defined on the factors. The
curvature for the product connection is explicitly calculated, and shown to be
independent of the choice of the twisting map and the module twisting map
used to define the product connection. As a consequence, we obtain that a
product of two flat connections is again a flat connection. We show that our
constructions also behaves well with respect to bimodule structures, namely
being the product of two bimodule connections again a bimodule connection. As
an application of our theory, all the product connections on the quantum plane
are computed.
\end{abstract}

\maketitle
\section*{Introduction}
\setcounter{equation}{0}

One of the main tools in classical differential geometry is the use
of the tangent bundle associated to a manifold. The role of the
algebra of functions on the manifold is taken by the sections of the
tangent bundle, namely, the vector fields. As a dual of the vector
fields space, the algebra of differential forms (endowed with the
exterior product) turns out to be an useful tool in the study of
global properties of the manifold, giving rise to invariants such as
the de Rham cohomology. A problem arises when trying to compare
vector fields and differential forms at different points of the
manifold, the solution to it being given by the concepts of
\dtext{(linear) connection} and \dtext{covariant derivative}, that
allow us to define the derivative of a curve on a point of orders higher
than one, hence giving us a way to speak about accelerations on a path. The
notion of connection also has another meanings in physics, like the existence of
an electromagnetic potential, which is equivalent to the existence of a
connection in a rank one trivial bundle with fixed trivialization.

Jean--Louis Koszul gave in \cite{Koszul60a} a powerful algebraic
generalization of differential geometry, in particular giving a
completely algebraic description of the notion of connection. These
notions were extended to a noncommutative framework by Alain Connes
in \cite{Connes86a}, what meant the dawn of noncommutative
differential geometry. Much research has been done about the theory
of connections in this context. On the one hand, Joachim Cuntz and
Daniel Quillen, in their seminal paper \cite{Cuntz95a} started the
theory of \dtext{quasi-free} algebras (also named \dtext{formally
smooth} by Maxim Kontsevich or \dtext{qurves} by Lieven Le Bruyn),
opening the way to an approach to noncommutative (algebraic)
geometry (also dubbed \dtext{nongeometry} to avoid confusions with
Michael Artin and Michel Van den Bergh's style of noncommutative
algebraic geometry). These formally smooth algebras are
characterized by the projectiveness (as a bimodule) of the first
order universal differential calculus, or equivalently as those
algebras that admit a universal linear connection. On the other
hand, in Connes' style of noncommutative geometry, the study of the
general theory of connections leads to the definition of the
Yang--Mills action, which turns out to be nothing but the usual
gauge action when we specialize it to the commutative case (cf.
\cite{Connes86a}, \cite{Landi97a}, \cite{Gracia01a} and references
therein).

In this paper, we deal with the problem of building up products of those
connection operators. Basically, there are two different notions of
``\emph{product connection}'' that one might want to build. Firstly, one might
want to consider two different bundles over a manifold, each of them endowed
with a connection, and then try to build a product connection on the
(fibre) product bundle. A noncommutative version of this construction was
given by Michel Dubois--Violette and John  Madore in \cite{DuboisViolette99a},
\cite{Madore95a}. Further steps on this direction, including its
relations with the realization of vector fields as Cartan pairs as proposed by
Andrzej Borowiec in \cite{Borowiec96a}, have been given by Edwin Beggs in
\cite{BeggsUNa}. The other possible notion of product connection, and the one
with which we want to deal, refers to the consideration of the cartesian product
of two given manifolds, and the building of a connection of the bundle
associated to this product manifold.

Traditionally, when taking the passage from classical geometry to
(noncommutative) algebra, the product space is associated with some kind of
tensor product (the algebraic tensor product in the case of algebraic
varieties, the topological tensor one when dealing with topological manifolds).
In \cite{Cap95a}, Andreas Cap, Herman Schichl and Ji\v{r}í Van\v{z}ura pointed
out the limitations of
these approach and proposed a definition of ``noncommutative cartesian
product'' of spaces by means of the so--called \dtext{twisted tensor product}
of the algebras. A twisted tensor product is a particular case of the notion of
\dtext{distributive law} given by Jon Beck in \cite{Beck69a}, and may be
regarded as a sort of local version of a braiding in a braided monoidal
category. Some further insights on the interpretation of this algebraic
construction from a geometrical point of view has been given by Alfons Van
Daele and S. Van Keer in \cite{VanDaele94a}, Stanislaw Woronowicz in
\cite{Woronowicz96a},
and the author and some collaborators in \cite{JaraUNa} (cf. also
\cite{LopezUNa} for an interpretation of twisted tensor products from a
deformation theory point of view). Following these ideas, we will show how to
build a product connection on a twisted tensor product of two algebras.

In Section \ref{preliminaries} we recall the notion of a (right) connection on
an algebra, given as an operator $\nabla:E\to E\ot_A \O^1 A$, where $E$ is a
(right) $A$--module and $\O^1 A$ a first order differential calculus over $A$,
motivating our choice of the differential calculus and the modules on the
building of the product connection, for which an explicit formula is given in
the simplest case of the usual tensor product. Then, we recall the main notions
we need about twisted tensor products, defined by means of twisting maps
$R:B\ot A\to A\ot B$.

In Section \ref{construction} we give the basic definition of our object of
study, and prove that our definition actually yields a connection in the
product space, and that this connection boils down nicely to the classical
product connection in the commutative case.

To a connection we can always associate its \dtext{curvature operator},
obtained by squaring the extension of the connection to the whole differential
calculus. The curvature operator leads to the definition of \dtext{flat
connections} as those having 0 curvature. Flat connections have been used by
Philippe Nuss in \cite{Nuss97a} in relation with noncommutative descent theory,
and also by Edwin Beggs and Tomasz Brzezinski in \cite{Beggs05a}, where they are
interpreted as the
differential of a certain complex in order to build a noncommutative de Rham
cohomology with coefficients. We deal with the problem of describing the
curvature of our product connection in Section \ref{curvature}, stating our
main theorem, that gives us an explicit formula to compute the curvature for
the product connection in terms of the curvatures of the factors:
    \begin{equation*}
        \theta (e\ot b, a\ot f) = i_E(\theta^E(e))\cd b + a\cd i_F(\theta^F(f)).
    \end{equation*}
The most striking consequence of this theorem is the fact that the curvature
does not depend neither on the twisting map $R$ nor on the module twisting map
that we use to get the module structure, suggesting that the curvature remains
invariant under all the deformations obtained by means of a twisted tensor
product. As an immediate corollary, we have that the product of two flat
connections is again a flat connection.

In Section \ref{bimodules} we consider bimodule connections (in the sense
 introduced by Jihad Mourad in \cite{Mourad95a}) instead of one sided
connections,
and we find necessary and sufficient conditions for the product of two bimodule
connections to be a bimodule connection. We conclude, in Section \ref{examples},
by illustrating our theory giving a complete description of all the product
connections on the quantum plane $k_q[x,y]$.

\section{Preliminaries}\label{preliminaries}
\setcounter{equation}{0}
\subsection{Connections on algebras}\hfill

Let $A$ be an associative, unital algebra over a field $k$, and
$\O A =\bigoplus_{p\geq 0}\O^p A$ a differential calculus over
$A$, that is, a differential graded algebra generated, as a differential graded algebra, by $\O^0 A\cong
A$, with differential $d=d_A$. Let $E$ be a (right) $A$--module; a
\dtext{(right) connection} on $E$ is a linear mapping
\[ \nabla: E \lto E\otimes_{A}\O^1 A
\]
satisfying the (right) \dtext{Leibniz rule}:
\begin{equation}
\nabla(s\cdot a)=(\nabla s)\cdot a + s\otimes da\quad \forall\,s\in
E, a\in A.
\end{equation}
Under these conditions, the mapping $\nabla$ can be extended in a
unique way to an operator
\[ \nabla: E\otimes_A \O A\lto E\otimes_A\O A
\]
of degree 1, by setting
\begin{equation}\label{connectionextension}
\nabla(s\otimes \omega)=\nabla{s}\otimes \omega + (-1)^{p}s\otimes
d\omega\quad\forall\,s\in E,\omega\in\O^p A,
\end{equation}
where we are using the identification $(E\ot_A \O^1 A)\ot_A \O^n
A\cong E\ot_A \O^{n+1}A$. Regarding $E\ot_A \O A$ as a right
$\O A$--module, we find that the following graded Leibniz rule is satisfied:
\begin{equation}
\nabla(\sigma\omega)=(\nabla\sigma)\omega + (-1)^p\sigma
d\omega\quad\forall\,\sigma\in E\ot_A \O^p A, \omega\in \O A.
\end{equation}
There are analogous concepts for left modules.

Usually, we will be interested on working with the universal
differential calculus over an algebra $A$. Connections over the
universal differential calculus will be called \dtext{universal
connections}. It is a well known fact (cf. \cite[Corollary 8.2]{Cuntz95a}) that a
right $A$--module admits a universal connection if, and only if, it
is projective over $A$.

Whenever $A$ is a commutative algebra, the tensor product $E\ot_A F$
of two $A$--modules $E$ and $F$ is again an $A$--module. If $E$ and
$F$ carry respective connections $\nabla^E$ and $\nabla^F$, we may
build the \dtext{tensor product connection} on $E\ot_A F$ by
defining
\begin{equation}
\nabla^{E\ot_A F}:=\nabla^E\ot F + E\ot \nabla^F.
\end{equation}
A possible generalization of this construction was given by
Dubois--Violette and Madore in \cite{DuboisViolette99a}, \cite{Madore95a}. If $E$ and $F$ are
$A$--bimodules equipped with right connections $\nabla^E$ and
$\nabla^F$, and such that there exists a linear mapping
\[ \sigma: \O^1 A\ot_A F \lto F\ot_A \O^1 A
\]
satisfying that
\begin{equation}
\nabla^F(am)=a\nabla^F(m) + \sigma(da\ot_A m)\quad\forall\,a\in
A,m\in F,
\end{equation}
then we may define
\[
    \nabla^{E\ot_A F}:E\ot_A F  \lto  E\ot_A F \ot \O^1 A
\]
by setting
\begin{equation}
    \nabla^{E\ot_A F}:= (E\ot \sigma)\circ(\nabla^E\ot F) + E\ot
    \nabla^F,
\end{equation}
and this $\nabla^{E\ot_A F}$ is a right connection on $E\ot_A F$.

Our aim is to define a different kind of ``\emph{product
connection}'' with a more geometrical flavour. Namely, consider that
our algebras $A=C^\infty(M)$ and $B=C^\infty(N)$ represent the
algebras of functions over certain manifolds $M$ and $N$, and that
$E=\cam(M)$ and $F=\cam(N)$ are the modules of vector fields on the
manifolds. The algebra associated to the cartesian product of the
manifolds is $C^\infty(M\times N)\cong C^\infty(M)\ot C^\infty (N)$
(more precisely, a suitable completion of the latest). For the
modules of vector fields and differential $1$--forms, we have that
\begin{gather*}
\cam(M\times N)\cong \cam(M)\ot C^\infty(N) \oplus C^\infty(M)\ot
\cam(N),\\
\O^1 (C^\infty(M)\ot C^\infty(N))\cong \O^1(C^\infty(M))\ot
C^\infty(N) \oplus C^\infty(M)\ot \O^1(C^\infty(N)),
\end{gather*}
hence, a ``\emph{product connection}'' of two connections defined on
$E$ and $F$ should be defined as a linear mapping
\[ \nabla : E\ot B \oplus A\ot F \lto  (E\ot B \oplus A\ot
F)\ot_{A\ot B} (\O^1 A\ot B \oplus A \ot \O^1B)
\]

Firstly, realize that if $E$ is a right (resp. left) $A$--module,
and $F$ is a right (resp. left) $B$--module, then $E\ot B\oplus A\ot
F$ is a right $(A\ot B)$--module, with actions
\begin{gather*}
(e\ot b,a\ot f)\cd(\alpha\ot \beta):=(e\alpha\ot b\beta,a\alpha\ot
f\beta)\\
\text{(resp.}\quad (\alpha\ot \beta)\cd (e\ot b,a\ot f):=(\alpha
e\ot \beta b,\alpha a\ot \beta f)\ \text{)}
\end{gather*}

For simplicity, we will only work with right connections. Left
connections admit a similar treatment.

\subsection{Product Connection}\hfill

Suppose then that $E$ is a right $A$--module endowed with a (right)
connection $\nabla^E$, and that $F$ is a right $B$--module endowed
with a (right) connection $\nabla^F$. Let us consider the mappings
\begin{eqnarray*}
\nabla_1:E\ot B & \lto & (E\ot B\op A\ot F)\ot_{A\ot B}(\O^1 A\ot
B\op A\ot \O^1 B), \\
\nabla_2:A\ot F & \lto & (E\ot B\op A\ot F)\ot_{A\ot B}(\O^1 A\ot
B\op A\ot \O^1 B)
\end{eqnarray*}
respectively given by
\begin{eqnarray*}
\nabla_1 &:=& (E\ot \tau \ot u_B)\circ (\nabla^E\ot B) + (E\ot
u_A\ot u_B\ot \O^1 B)\circ (E\ot d_B),\ \text{and}\\
\nabla_2 &:=& (A\ot F\ot u_A\ot \O^1B)\circ (A\ot \nabla^F) +
(u_A\ot \tau\ot u_B)\circ (d_A\ot F),
\end{eqnarray*}
where $\tau$ represent classical flips. If we use the shorthand
notation $\nabla^E(e)=e_i\ot d_A a_i$, where the summation symbol is
omitted, the Leibniz rule for $\nabla^E$ is written as
\begin{equation}\label{Leibiz:E}
\nabla^E(e\alpha)=e_i\ot (d_A a_i)\alpha + e\ot d\alpha,
\end{equation}
and we have that
\begin{eqnarray*}
    \nabla_1((e\ot b)\cd(\alpha\ot \beta)) &=& \nabla_1(e\alpha\ot
    b\beta) = \\
    &=& e_i\ot b\beta \ot_{A\ot B} (d_A a_i)\alpha \ot 1 + e\ot b\beta
    \ot_{A\ot B} d\alpha \ot 1 + \\
    & & +\  e\alpha \ot 1 \ot_{A\ot B} 1 \ot d_B(b\beta)= \\
    &=& e_i\ot b \ot_{A\ot B} (d_A a_i)\alpha \ot \beta + e\ot b
    \ot_{A\ot B} d\alpha \ot \beta + \\
    & & +\  e \ot 1 \ot_{A\ot B} \alpha \ot d_B(b)\beta + e \ot 1 \ot_{A\ot B} \alpha \ot b d_B \beta = \\
    &=& \left(e_i\ot b \ot_{A\ot B} da_i\ot 1 + e\ot 1\ot_{A\ot B} 1\ot
    db)\cd (\alpha\ot \beta \right) + \\
    && +\ e\ot b \ot_{A\ot B} d_A\alpha\ot \beta + e\ot b \ot_{A\ot B} \alpha\ot
    d_B\beta = \\
    &=& \nabla_1(e\ot b)\cd (\alpha\ot \beta) + (e\ot b) \ot_{A\ot
    B} d(\alpha\ot \beta).
\end{eqnarray*}
A similar computation shows that
\begin{eqnarray*}
    \nabla_2((a\ot f)\cd(\alpha\ot \beta)) &=& \nabla_2(a\ot f)\cd (\alpha\ot \beta) + (a\ot f) \ot_{A\ot
    B} d(\alpha\ot \beta).
\end{eqnarray*}
Adding up these two equalities, we conclude that the map
\begin{eqnarray*}
\nabla : E\ot B \oplus A\ot F & \lto & (E\ot B \oplus A\ot
F)\ot_{A\ot B} (\O^1 A\ot B \oplus A \ot \O^1B)\\
(e\ot b,a\ot f) & \lmto & \nabla_1(e\ot b) + \nabla_2(a\ot f)
\end{eqnarray*}
verifies that
\[ \nabla((e\ot b,a\ot f)\cd (\alpha\ot \beta))= \nabla(e\ot b,a\ot f)\cd (\alpha\ot
\beta) + (e\ot b,a\ot f) \ot_{A\ot B} d(\alpha \ot \beta),
\]
and henceforth, $\nabla$ is a (right) connection on the module $E\ot
B\op A\ot F$. We shall call this map the \dtext{(classical) product
connection} of $\nabla^E$ and $\nabla^F$.

\subsection{Twisted tensor products}\hfill

Let $k$ be a field, used as a base field throughout. We denote
$\otimes _k$ by $\otimes $, the identity $id_V$ of an object $V$
simply by $V$, and by $\tau :V\otimes W \rightarrow W\otimes V$,
$\tau (v\otimes w)=w\otimes v$, the usual flip. All algebras are
assumed to be associative unital $k$-algebras; the multiplication
and unit of an algebra $D$ are denoted by $\mu _D:D\otimes
D\rightarrow D$ and respectively $u_D:k\rightarrow D$ (or simply by
$\mu $ and $u$ if there is no danger of confusion).

 We recall the twisted tensor product of algebras from
\cite{Tambara90a}, \cite{VanDaele94a}, \cite{Cap95a}. If $A$ and $B$ are two
algebras, a linear map $R:B\ot A\rightarrow A\ot B$ is called a
\dtext{twisting map} if it satisfies the conditions
\begin{eqnarray}
&&R(b\ot 1)=1\ot b,\;\;\;R(1\ot a)=a\ot 1, \;\;\;\forall \; a\in A,
\;b\in B, \label{tw1}\\
&&R\circ (B\ot \mu _A)=(\mu _A\ot B)\circ (A\ot R)\circ (R\ot A),
\label{tw2}\\
&&R\circ (\mu _B\ot A)=(A\ot \mu _B)\circ (R\ot B)\circ (B\ot R).
\label{tw3}
\end{eqnarray}
If we denote by $R(b\ot a)=a_R\ot b_R$, for $a\in A$, $b\in B$, then
(\ref{tw2}) and (\ref{tw3}) may be written as:
\begin{eqnarray}
&&(aa')_R\ot b_R=a_Ra'_r\ot (b_R)_r, \label{tw4} \\
&&a_R\ot (bb')_R=(a_R)_r\ot b_rb'_R, \label{tw5}
\end{eqnarray}
for all $a, a'\in A$ and $b, b'\in B$, where $r$ is another copy of $R$.
If we define a multiplication on $A\ot B$, by $\mu _R=(\mu _A\otimes
\mu _B)\circ (A\otimes R\otimes B)$, that is
\begin{eqnarray}
&&(a\ot b)(a'\ot b')=aa'_R\ot b_Rb', \label{multtw}
\end{eqnarray}
then this multiplication is associative and $1\ot 1$ is the unit.
This algebra structure is denoted by $A\ot _RB$ and is called the
\dtext{twisted tensor product} of $A$ and $B$. This construction
works also if $A$ and $B$ are algebras in an arbitrary monoidal
category.

 If $A\otimes _{R_1} B$, $B\otimes _{R_2} C$ and $A\otimes
_{R_3} C$ are twisted tensor products of algebras, the twisting maps
$R_1$, $R_2$, $R_3$ are called {\it compatible} if they satisfy
\begin{eqnarray*}
&&(A\otimes R_2)\circ (R_3\otimes B)\circ (C\otimes R_1)=
(R_1\otimes C)\circ (B\otimes R_3)\circ (R_2\otimes A),
\end{eqnarray*}
see \cite{JaraUNa}. If this is the case, the maps $T_1:C\otimes
(A\otimes _{R_1} B)\rightarrow (A\otimes _{R_1} B)\otimes C$ and
$T_2:(B\otimes _{R_2} C)\otimes A\rightarrow A\otimes (B\otimes
_{R_2} C)$ given by $T_1:=(A\otimes R_2)\circ (R_3\otimes B)$ and
$T_2:=(R_1\otimes C)\circ (B\otimes R_3)$ are also twisting maps and
$A\otimes _{T_2} (B\otimes _{R_2} C)\equiv (A\otimes _{R_1} B)
\otimes _{T_1} C$; this algebra is denoted by $A\otimes _{R_1} B
\otimes _{R_2} C$. This construction may be iterated to an arbitrary
number of factors, see \cite{JaraUNa} for complete detail.

When we have a left $A$--module $M$, a left $B$--module $N$, a
twisting map $R:B\otimes A\to A\otimes B$ and a linear map
$\tau_{M,B}:B\otimes M \to M\otimes B$ such that
\begin{eqnarray}\label{moduletwistcondition1}
        \tau_{M,B}\circ(\mu_B\otimes M) & = & (M\otimes \mu_B)\circ (\tau_{M,B}\otimes B)\circ (B\otimes
        \tau_{M,B}), \\ \label{moduletwistcondition2}
        \tau_{M,B}\circ(B\otimes \l_M) & = & (\l_M \otimes B)\circ (A\otimes \tau_{M,B})\circ (R\otimes
        M),
\end{eqnarray}
then the map $\l_{\tau_{M,B}}:(A\otimes_R B)\otimes (M\otimes N)\to
M\otimes N$ defined by $\l_{\tau_{M,B}}:=(\l_M\otimes \l_N)\circ
(A\otimes \tau_{M,B}\otimes N)$ yields a left $(A\otimes_R
B)$--module structure on $M\otimes N$, which furthermore is
compatible with the inclusion of $A$. In this case, we say that
$\tau_{M,B}$ is a \dtext{(left) module twisting map}. Unlike what happens for algebra twisting
maps, usually is not enough to have a left $(A\otimes_R B)$--module
structure on $M\otimes N$ in order to recover a module twisting map.
Some sufficient conditions for this to happen are, for instance, requiring
that $M$ is projective and $N$ is faithful (cf. \cite[Theorem 3.8]{Cap95a}).

Similarly, if we have a twisting map $R:B\otimes A\to A\otimes B$, a right $A$--module $M$ and a
right $B$--module $N$, a linear map $\tau_{N,A}:N\otimes A\to A\otimes N$ such that
\begin{eqnarray}\label{rightmoduletwist1}
    \tau_{N,A}\circ (N\otimes \mu_A) &=&(\mu_A\otimes N)\circ (A\otimes
    \tau_{N,A})\circ (\tau_{N,A}\otimes A)\\
    \label{rightmoduletwist2}
    \tau_{N,A}\circ (\rho_B\otimes A) &=& (A\otimes \rho_B)\circ
    (\tau_{N,A}\otimes B)\circ (N\otimes R).
\end{eqnarray}
then the map $\rho_{\tau_{N,A}}:=(\rho_A\otimes \rho_B)\circ (M\otimes \tau_{N,A}\otimes B)$, yields a right $(A\ot_R B$--module action on $M\ot N$. In this case, we call $\tau_{N,A}$ a \dtext{(right) module twisting map}

Twisting maps also have a nice behaviour with respect to (universal) differential calculi. More concretely, we have the following result (cf. \cite{Cap95a}):

\begin{theorem}\label{twistdifferentialforms}
Let $A$, $B$ be two algebras. Then any twisting map $R:B\otimes A\to
A\otimes B$ extends to a unique twisting map $\tilde{R}:\O B\otimes
\O A \to \O A\otimes \O B$ which satisfies the conditions
\begin{eqnarray}\label{twisteddiff1}
\tilde{R}\circ(d_B\otimes \O A)&=&(\eps_A \otimes d_B)\circ
\tilde{R}, \\
\label{twisteddiff2} \tilde{R}\circ (\O B\otimes d_A)&=&(d_A \otimes
\eps_B)\circ \tilde{R},
\end{eqnarray}
where $d_A$ and $d_B$ denote the differentials on the algebras of
universal differential forms $\O A$ and $\O B$,
and $\eps_A$, $\eps_B$ stand for the gradings on $\O A$ and $\O B$,
respectively. Moreover, $\O A\otimes_{\tilde{R}}\O B$ is a graded
differential algebra with differential $d(\vphi\otimes
\omega):=d_A\vphi\otimes \omega + (-1)^{\abs{\vphi}}\vphi\otimes
d_B\omega$.
\end{theorem}

\section{Twisted tensor product connection}\label{construction}
\setcounter{equation}{0}

In the former section we introduced the definition of a connection
within the formalism of differential calculus over algebras, and
showed how to build the product connection for a tensor product of
two algebras, extending the definition of the classical product
connection in differential geometry. In \cite{JaraUNa}, we advocated
that in noncommutative geometry the cartesian product should not be
replaced at the algebraic level by the usual tensor product of
algebras, but by a deformation of it, known as the twisted tensor
product. In this section, we will show how to extend the definition
of the product connection to a twisted tensor product of two
algebras under suitable conditions.

Let $A$ and $B$ be algebras, $R:B\ot A\to A\ot B$ a twisting map,
$E$ a right $A$--module endowed with a right connection $\nabla^E$,
and $F$ a right $B$--module endowed with a right connection
$\nabla^F$. In \cite{Cap95a} it is shown that we can lift the
twisting map $R$ to a twisting map $\widetilde{R}:\O B\ot \O A\to \O
A\ot \O B$ on the graded differential algebras of (universal)
differential forms, and that the algebra
\[ \O A\ot_{\widetilde{R}} \O B = \bigoplus_{n\in \n}\left(\bigoplus_{p+q=n}\O^p A\ot \O^q B\right)
\]
is a differential calculus over $A\ot_{R} B$. For this differential
calculus, the module of $1$--forms can be identified as $\O^1 A \ot
B \oplus A\ot \O^1 B$, with the natural action induced by the
twisting map. As the situation is pretty much the same as in the
tensor product case, the natural way for defining a ``\emph{twisted
product}'' connection of $\nabla^E$ and $\nabla^F$ would be
considering a linear map
\[ \nabla: E\ot B \op A\ot F \lto \left(E\ot B \op A\ot
F\right)\ot_{A\ot_{\widetilde{R}}B} \left(\O^1 A\ot B \op A\ot \O^1
B\right).
\]
The first step on making this map becoming a connection is giving a
right $(A\ot_R B)$--module action on $E\ot B \op A\ot F$, which
means finding a right $(A\ot_R B)$--module structure on both $E\ot
B$ and $A\ot F$. For the first one we may just use the twisting map
and define:
\begin{equation} \label{eq:twistactionEB}
    (e\ot b)\cd (\alpha\ot \beta):= e\alpha_R \ot b_R\beta.
\end{equation}
For the second one, a sufficient way of giving a module structure is
finding a (right) module twisting map $\tau_{F,A}:F\ot A\to A\ot F$,
and then taking
\begin{equation} \label{eq:twistactionAF}
    (a\ot f)\cd(\alpha\ot \beta):=a\alpha_\tau\ot f_\tau\beta.
\end{equation}
The fact that the former definitions are indeed module actions
follows directly  from the fact that both $R$ and $\tau_{F,A}$ are
right module twisting maps (cf. \cite{Cap95a}, 3.12).

Following the lines given by the definition of the classical tensor
product connection, in order to build $\nabla$ we have to find
suitable maps $\nabla_1$ and $\nabla_2$. For the first one, it
suffices to define

\begin{gather*}
\nabla_1:E\ot B  \lto  (E\ot B\op A\ot F)\ot_{A\ot B}(\O^1 A\ot
B\op A\ot \O^1 B)\\
\nabla_1:=(E\ot u_B\ot \O^1 A\ot B)\circ (\nabla^E\ot B) + (E\ot
u_B\ot u_A\ot \O^1 B)\circ (E\ot d_B).
\end{gather*}
With this definition, when $R$ is the classical flip is a definition
trivially equivalent to the one given in the former section, and we
have that

\begin{eqnarray*}
    \nabla_1((e\ot b)\cd(\alpha\ot \beta)) &=& \nabla_1(a\alpha_R\ot
    b_R\beta) = \\
    &=& (E\ot u_B\ot \O^1 A\ot B)(\nabla^E(e\alpha_R)\ot b_R\beta)+\\
    && +(E\ot u_B\ot u_A\ot \O^1 B)(e\alpha_R\ot d(b_R\beta))\overset{1}{=}\\
    &\overset{1}{=}& e_i\ot 1\ot_{A\ot_R B} (d_Aa_i)\alpha_R\ot b_R\beta +\\
    && +\ e\ot 1 \ot_{A\ot_R B} d\alpha_R\ot b_R\beta
    + \\
    && +\ e\alpha_R\ot 1\ot_{A\ot_R B}1\ot (d_B b_R)\beta + \\
    && +  e\alpha_R\ot 1\ot_{A\ot_R B}1\ot b_R
    d_B\beta=\\
    &=& e_i\ot 1\ot_{A\ot_R B} (d_Aa_i)\alpha_R\ot b_R\beta +\\
    && +\ e\ot 1 \ot_{A\ot_R B} d\alpha_R\ot b_R\beta
    + \\
    &&+\ e\ot 1\ot_{A\ot_R B}\alpha_R\ot (d_B b_R)\beta + \\
    && + e\ot 1\ot_{A\ot_R B}\alpha_R\ot b_R
    d_B\beta \overset{2}{=}\\
    &\overset{2}{=}& (e_i\ot 1\ot_{A\ot_R B}d_Aa_i \ot b + \\
    && +\ e\ot 1\ot_{A\ot_R
B}1\ot
    b)\cd (\alpha\ot \beta)+\\
    &&+\ e\ot b \ot_{A\ot_R B} d_A\alpha\ot \beta +\\
    && +\ e\ot b \ot_{A\ot_R B}
\alpha\ot
    d_B\beta =\\
    &=& \nabla_1(e\ot b)\cd (\alpha\ot \beta) + e\ot b \ot_{A\ot_R
    B} d(\alpha\ot\beta),
\end{eqnarray*}
where in 1 we are using Leibniz's rules (for the connection
$\nabla^E$ and the differential $d_B$), in $2$ the definition of the
action \eqref{eq:twistactionEB} and the compatibility of the
twisting map with the differential, as mentioned in equations \eqref{twisteddiff1} and \eqref{twisteddiff2}.

The definition of $\nabla_2$ is more involved, and we are forced to
assume some extra conditions on the maps $R$ and $\tau_{F,A}$.
Namely, assume that $R$ is invertible, with inverse $S:A\ot B\to
B\ot A$, that $\tau_{F,A}$ is invertible with inverse
$\sigma_{A,F}:A\ot F\to F\ot A$, and such that the following
relation, ensuring the compatibility of the module
twisting map with the connection $\nabla^F$, is satisfied:
\begin{gather}
\label{eq:nabla2cond1} (A\ot \nabla^F)\circ
\tau_{F,A}=(\tau_{F,A}\ot\O^1 B)\circ(F\ot
\widetilde{R})\circ(\nabla^F\ot A).
\end{gather}
From this condition, that in Sweedler's like notation is written as
    \begin{equation}
        a_\tau\ot (f_\tau)_j\ot_B (db_\tau)_j =
        (a_{\widetilde{R}})_\tau\ot (f_j)_\tau\ot_B ((db_j)_{\widetilde{R}})_\tau,
    \end{equation}
the module twisting conditions \eqref{rightmoduletwist1} and \eqref{rightmoduletwist2} for
$\tau_{F,A}$,  and the twisting map
conditions \eqref{tw2} and \eqref{tw3} for $R$, we may
easily deduce the following equalities:
\begin{gather}
\label{eq:nabla2cond1bis} (\sigma_{A,F} \ot \O^1 B)\circ (A\ot\nabla^F) = (F\ot\widetilde{R})\circ(\nabla^F\ot A)\circ \sigma_{A,F},\\
\label{eq:nabla2cond2} (\mu_A\ot F)\circ(A\ot
\tau_{F,A})=\tau_{F,A}\circ(F\ot \mu_A)\circ (\sigma_{A,F}\ot A),\\
\label{eq:nabla2cond3} \sigma_{A,F}\circ(A\ot
\lambda_F)\circ(\tau_{F,A}\ot B) = (\lambda_F\ot A)\circ(F\ot S), \\
\label{eq:nabla2cond4} \sigma_{A,F}\circ (\mu_A\ot F) = (F\ot
\mu_A)\circ
(\sigma_{A,F}\ot A)\circ (A\ot \sigma_{A,F}),\\
\label{eq:nabla2cond5} \sigma_{A,F}\circ(B\ot \lambda_F)\circ(R\ot
F) = (\lambda_F\ot A)\circ (B\ot \sigma_{A,F})
\end{gather}
and define the map
\begin{gather*}
    \nabla_2: A\ot F  \lto  (E\ot B\op A\ot F)\ot_{A\ot B}(\O^1 A\ot
B\op A\ot \O^1 B) \\
    \nabla_2:= (A\ot F\ot u_B\ot \O^1 B)\circ(A\ot \nabla^F) +
    (u_A\ot F\ot d_A\ot u_B)\circ \sigma
\end{gather*}
then we have that
\begin{eqnarray*}
    \nabla_2((a\ot f)\cd(\alpha\ot\beta)) &=&
    \nabla_2(a\alpha_\tau\ot f_\tau\beta)= \\
    &=& (A\ot F\ot u_A\ot \O^1
    B)(a\alpha_\tau\ot\nabla^F(f_\tau\beta)) + \\
    && + \ 1\ot (f_\tau\beta)_\sigma\ot
    d_A((a\alpha_\tau)_\sigma)\ot 1 \overset{\eqref{eq:nabla2cond4}}{=}\\
    &\overset{\eqref{eq:nabla2cond4}}{=}& a \alpha_\tau\ot
    (f_\tau)_j\ot_{A\ot_R B} 1 \ot d_B(b_\tau)_j\beta + \\
    && +\ a\alpha_\tau\ot f_\tau\ot_{A\ot_R B} 1\ot d_B\beta + \\
    && +\ 1\ot (f_\tau\beta)_{\sigma\overline{\sigma}}\ot_{A\ot_R B}
    (d_Aa_{\overline{\sigma}})\alpha_{\tau\sigma}\ot 1 + \\
    &&+\ 1\ot (f_\tau\beta)_{\sigma\overline{\sigma}}\ot_{A\ot_R B}
    a_{\overline{\sigma}}d_A(\alpha_{\tau\sigma})\ot 1
    \overset{\eqref{eq:nabla2cond1}}{=}\\
    &\overset{\eqref{eq:nabla2cond1}}{=}& a (\alpha_{\widetilde{R}})_{\tau}\ot
    (f_j)_\tau \ot_{A\ot_R B} 1 \ot (d_Bb_j)_{\widetilde{R}}\beta + \\
    && +\ a \ot f\ot_{A\ot_R B} \alpha\ot d_B\beta + \\
    && +\ 1\ot (f_\tau\beta)_{\sigma\overline{\sigma}}\ot_{A\ot_R B}
    (d_Aa_{\overline{\sigma}})\alpha_{\tau\sigma}\ot 1 + \\
    &&+\ a \ot (f_\tau\beta)_{\sigma}\ot_{A\ot_R B}
   d_A(\alpha_{\tau\sigma})\ot 1 \overset{\eqref{eq:nabla2cond3}}{=}\\
   & \overset{\eqref{eq:nabla2cond3}}{=} & a \ot
    f_j \ot_{A\ot_R B} \alpha_{\widetilde{R}} \ot (d_Bb_j)_{\widetilde{R}}\beta
+ \\
    && +\ a \ot f\ot_{A\ot_R B} \alpha\ot d_B\beta + \\
    && +\ 1\ot (f\beta_S)_{\sigma}\ot_{A\ot_R B}
    (d_Aa_{\sigma})\alpha_{S}\ot 1 + \\
    &&+\ a \ot f\beta_S \ot_{A\ot_R B}
   d_A(\alpha_{S})\ot 1 \overset{\eqref{eq:nabla2cond2}}{=}\\
   &\overset{\eqref{eq:nabla2cond2}}{=}& (a \ot
    f_j \ot_{A\ot_R B} 1\ot (d_Bb_j))\cd(\alpha\ot\beta) + \\
    && +\ a \ot f\ot_{A\ot_R B} \alpha\ot d_B\beta + \\
    && +\ 1\ot f_\sigma\beta_{S\overline{S}}\ot_{A\ot_R B}
    d_A(a_{\sigma\overline{S}})\alpha_{S}\ot 1 + \\
    && +\ a \ot f \ot_{A\ot_R B} d_A\alpha\ot \beta =\\
   &=& (a \ot
    f_j \ot_{A\ot_R B} 1\ot (d_Bb_j))\cd(\alpha\ot\beta) + \\
    && +\ 1\ot f_\sigma\ot_{A\ot_R B}
    d_A(a_{\sigma})\alpha\ot \beta +     \\
    && +\ a \ot f\ot_{A\ot_R B} \alpha\ot d_B\beta + \\
    && +\ a \ot f \ot_{A\ot_R B}
   d_A\alpha\ot \beta =\\
   &=& \nabla_2(a\ot f)\cd(\alpha\ot \beta) + a\ot f\ot_{A\ot_R
   B}d(\alpha\ot \beta).
\end{eqnarray*}
Henceforth, the mapping
\[
    \nabla: E\ot B \op A\ot F \lto \left(E\ot B \op A\ot F\right)\ot_{A\ot_{\widetilde{R}}B} \left(\O^1 A\ot B \op A\ot \O^1 B\right)
\]
defined as
\begin{equation}\label{productconnection}
\nabla(e\ot b,a\ot f) := \nabla_1(e\ot b) + \nabla_2(a\ot f)
\end{equation}
is a (right) connection on the module $E\ot
B\oplus A\ot F$. We will call this connection the \dtext{(twisted) product connection of $\nabla^E$ and $\nabla^F$}.

\section{Curvature on product connections}\label{curvature}
\setcounter{equation}{0}

In this section our aim is to study the curvature for the formerly defined product connections. If we have a connection $\nabla:E \to E \ot_A \O^1 A$, we will also denote by $\nabla:E\ot_A \O A\to E\ot_A \O A$ the extension given by \eqref{connectionextension}, ocassionally denoting by $\nabla^{[n]}:E\ot_A \O^n A\to E\ot_A \O^{n+1}A$ its restriction to $E$-valued $n$--form. The \dtext{curvature} of the connection $\nabla$ is defined to be the operator $\theta:=\nabla^{[1]}\circ \nabla^{[0]}:E\to E\ot_A\O^2 A$. It is well known (cf. for instance \cite[Sect. 7.2]{Landi97a}) that the map $\theta$ is right $A$--linear. A connection $\nabla$ is said to be a \dtext{flat connection} whenever the associated curvature map is equal to 0. As curvature map may be extended to a (right) $\O A$--linear map $\theta:E\ot_A \O A \to E\ot_A \O A$ of degree 2 given at degree $n$ by $\theta^{[n]}:=\nabla^{[n+1]}\circ \nabla^{[n]}$, and it is easily checked that $\theta^{[n]}= (E\ot\mu_{\O A}) \circ (\O^n A\ot \theta)$, (cf. \cite[Prop 2.3]{Beggs05a}), we have that a flat connection  can be used for building a noncommutative de Rham cohomology with a nontrivial coefficient bundle.

Let then $A$ and $B$ algebras, $R:B\ot A\to A\ot B$ a twisting map,
$E$ a right $A$--module endowed with a right connection $\nabla^E$,
and $F$ a right $B$--module endowed with a right connection
$\nabla^F$ such that we can build the product connection $\nabla$ as
in the former section, let also $\nabla=(\nabla^{[n]})$ denote the
extension of $\nabla$ to $(E\ot B\oplus A\ot F)\ot_{A\ot_R B}(\O
A\ot_{\widetilde{R}} \O B)$. For $e\in E$, let us denote
$\nabla^E(e)=e_i\ot_A d_A a_i$, and $\nabla^E(e_i):=e_{ij}\ot_A
d_Aa_{ij}$, where summation symbols are omitted. In the same spirit,
for $f\in F$, we will denote $\nabla^F(f)=f_k\ot_B d_B b_k$, and
$\nabla^F(f_k):=f_{kl}\ot_B d_B b_{kl}$. With this notation, the
respective curvatures are written as $\theta^E(e)=e_{ij}\ot_A
d_Aa_{ij}d_A{a_i}$, $\theta^F(f)=f_{kl}\ot_B d_B b_{kl} d_B{b_k}$.
We will also denote by $i_E$ and $i_F$ the canonical inclusions (as
vector spaces) of $E\ot_A \O^2 A$ and $F\ot_B \O^2 B$ into $(E\ot
B\oplus A\ot F)\ot_{A\ot_{R}B}(\O A\ot_{\widetilde{R}}\O B)^2$. For
a generic element $(e\ot b,a\ot f)\in (E\ot B\oplus A\ot F)$, using
the definition of the product connection we have that
\begin{eqnarray*}
    \nabla(e\ot b,a\ot f) &=& e_i\ot 1\ot_{A\ot_{R}B} d_Aa_i\ot b + e\ot 1\ot_{A\ot_{R}B} 1\ot d_Bb + \\
    && + 1\ot f_\sigma \ot_{A\ot_{R}B} d_A(a_\sigma)\ot 1 + a\ot f_k\ot_{A\ot_{R}B} 1\ot d_Bb_k
\end{eqnarray*}
Applying $\nabla^{[1]}$ to each of these four term we obtain:
\begin{eqnarray*}
        \nabla^{[1]}(e_i\ot 1\ot_{A\ot_{R}B} d_Aa_i\ot b) & = & \nabla(e_i\ot
1)\cd (d_Aa_i\ot b) +\\
        && +\ (e_i\ot 1)\ot_{A\ot_{R}B} d(da_i\ot b) \overset{1}{=}\\
        & \overset{1}{=} &(e_ij\ot 1 \ot_{A\ot_{R}B} d_Aa_{ij}\ot 1)\cd
(d_Aa_i\ot b) \\
        && - e_i\ot 1\ot_{A\ot_{R}B} d_Aa_i\ot d_Bb = \\
        &= & e_{ij}\ot 1 \ot_{A\ot_{R}B} d_Aa_{ij}d_A a_i\ot b \\
        && - e_i\ot 1\ot_{A\ot_{R}B} d_Aa_i\ot d_Bb = \\
        & = & i_E(\theta^E(e))\cd b - e_i\ot 1\ot_{A\ot_{R}B} d_Aa_i\ot d_Bb,
\end{eqnarray*}
\begin{eqnarray*}
    \nabla^{[1]}(e\ot 1\ot_{A\ot_{R}B} 1\ot d_Bb) & = & \nabla(e\ot 1)\cd 1\ot d_Bb + (e\ot 1)\ot_{A\ot_{R}B} d(1\ot d_Bb) = \\
    & = & e_i\ot 1\ot_{A\ot_{R}B} d_Aa_i\ot d_Bb,
\end{eqnarray*}
\begin{eqnarray*}
    \nabla^{[1]}(1\ot f_\sigma \ot_{A\ot_{R}B} d_A(a_\sigma)\ot 1) & = & \nabla(1\ot f_\sigma)\cd (d_A(a_\sigma)\ot 1) + \\
    && + (1\ot f_\sigma)\ot_{A\ot_{R}B} d(d_A(a_\sigma)\ot 1) = \\
    &=& (1\ot(f_\sigma)_k\ot_{A\ot_{R}B} 1\ot d_B(b_\sigma)_k)\cd (d_A (a_\sigma)\ot 1) = \\
    & = & 1\ot(f_\sigma)_k\ot_{A\ot_{R}B} (d_A(a_\sigma))_{\widetilde{R}}\ot (d_B(b_\sigma)_k)_{\widetilde{R}} \overset{2}{=}\\
    & \overset{2}{=} & -1\ot(f_\sigma)_k\ot_{A\ot_{R}B} d_A(a_\sigma)_{\widetilde{R}}\ot (d_B(b_\sigma)_k)_{\widetilde{R}}\overset{\eqref{eq:nabla2cond1bis}}{=} \\
    &\overset{\eqref{eq:nabla2cond1bis}}{=}& -1\ot(f_k)_\sigma \ot_{A\ot_{R}B} d_A(a_\sigma)\ot d_Bb_k,
\end{eqnarray*}
\begin{eqnarray*}
    \nabla^{[1]}(a\ot f_k\ot_{A\ot_{R}B} 1\ot d_Bb_k) & = & \nabla(a\ot
f_k)\cd(1\ot d_Bb_k) + \\
    && +\ a\ot f_k\ot_{A\ot_{R}B} d(1\ot d_Bb_k) = \\
    &=& (a\ot f_{kl} \ot_{A\ot_{R}B} 1\ot d_B b_{kl})\cd(1\ot d_Bb_k) + \\
    && + (1\ot (f_k)_\sigma \ot_{A\ot_{R}B} d_A a_\sigma\ot 1)\cd(1\ot d_Bb_k) =\\
    &=& a\ot f_{kl} \ot_{A\ot_{R}B} 1\ot d_B b_{kl}d_Bb_k + \\
    && + 1\ot (f_k)_\sigma \ot_{A\ot_{R}B} d_A a_\sigma \ot d_Bb_k = \\
    &=& a\cd i_F(\theta^F(f)) + 1\ot (f_k)_\sigma \ot_{A\ot_{R}B} d_A a_\sigma \ot d_Bb_k.
\end{eqnarray*}
where in $1$ we are using the definitions of $\nabla$ and the differential $d$, in $2$ the compatibility of $\widetilde{R}$ with $d_A$. Adding up these four equalities we obtain the following result:

\begin{theorem}
    The curvature of the product connection is given by
\begin{equation}\label{curvatureformula}
    \theta (e\ot b, a\ot f) = i_E(\theta^E(e))\cd b + a\cd i_F(\theta^F(f)).
\end{equation}
\end{theorem}

An interesting remark at the sight of the former result is that the product
curvature does not depend neither on the twisting map $R$ nor on the module
twisting map $\tau_{F,A}$, but only on the curvatures of the factors.
As an immediate consequence of Equation \eqref{curvatureformula} we obtain the following result:

\begin{corollary}
    The product connection of two flat connections is a flat connection.
\end{corollary}

Henceforth, one might ask the question of describing the de Rham
cohomology with coefficients in the sense of Beggs and Brzezinski
(ref. \cite{Beggs05a}) for the (twisted) product connection of two
flat connections. We will leave this problem for future works. It is
also worth noticing that formula \eqref{curvatureformula} drops down
in the commutative case to the classical formula for the curvature
on a product manifold.

\section{Bimodule connections}\label{bimodules}
\setcounter{equation}{0}

For many purposes, only considering right (or left) modules is not enough. On the one hand, if we want to apply our theory to $\ast$--algebras, then sooner or later we will be bond to deal with $\ast$--modules and hermitian modules, but since the involution reverses the order of the products, these notions only make sense when we consider bimodules. On the other hand, there is a special kind of connections, known as \dtext{linear connections}, obtained when we take $E=\O^1 A$. Since $\O^1 A$ is a bimodule in a natural way, there is no reason to neglect one of its structures restraining ourselves to look at it just as a one-sided module. Reasons for extending the notion of connection to bimodules have been largely discussed at \cite{Mourad95a}, \cite{DuboisViolette99a} and references therein.

Different approaches for dealing with this problem have been tried. The first one, described by Cuntz and Quillen in \cite{Cuntz95a}, consists on considering a couple $(\nabla^l,\nabla^r)$ where $\nabla^l$ is a left connection which is also a right $A$--module morphism, and $\nabla^r$ a right connection which is also a left $A$--module morphism. As it was pointed out in \cite{Dabrowski96a}, this approach, though rising a very interesting algebraic theory, is not well suited for our geometrical point of view, since it doesn't behave as expected when restricted to the commutative case. A different approach was introduced by Mourad in \cite{Mourad95a} for the particular case of linear connections and later generalized to arbitrary bimodules by Dubois-Violette and Masson in \cite{DuboisViolette96a} (see also \cite[Chapter 10]{DuboisViolette99a}). Their approach goes as follow: let $E$ be an $A$--bimodule; a \dtext{(right) bimodule connection} on $E$ is a right connection $\nabla:E\to E\ot_A\O^1 A$ together with a bimodule homomorphism $\sigma:\O^1 A\ot_A E \to E\ot_A \O^1 A$ such that
\begin{equation}
    \nabla(ma)=a\nabla(m) + \sigma(d_A(a)\ot_A m)\quad \text{for any $a\in A$, $m\in E$}.
\end{equation}
Giving a right bimodule connection in the above sense is equivalent to give a pair $(\nabla^L,\nabla^R)$ consisting in a left connection $\nabla^L$ and a right connection $\nabla^R$ that are $\sigma$--compatible, meaning that
\begin{equation}\label{sigmacompatible}
    \nabla^R = \sigma\circ \nabla^L.
\end{equation}

\begin{rem}
     A weaker definition of $\sigma$--compatibility, namely requiring that equation \eqref{sigmacompatible} holds only in the center $Z(E):=\{m\in E:\ am=ma\ \forall a\in A\}$ of $E$ rather than in the whole bimodule, has also been studied in \cite{Dabrowski96a}.
\end{rem}

So, assume that we have $E$ bimodule over $A$, $\nabla^E$ a bimodule connection on $E$ with respect to the morphism $\vphi:\O^1A\ot_A E\to E\ot_A\O^1 A$, and $F$ a bimodule over $B$ endowed with $\nabla^F$ a bimodule connection with respect to the bimodule morphism $\psi:\O^1 B\ot_B F\to F\ot_B\O^1 B$. As before, let $R:B\ot A\to A\ot B$ an invertible twisting map with inverse $S$, and assume also that we have a right module twisting maps $\tau_{F,A}:F\ot A\to A\ot F$ satisfying condition \eqref{eq:nabla2cond1} and a left module twisting map $\tau_{B,E}:B\ot E\to E\ot B$ satisfying condition
    \begin{equation}\label{eq:nabla1cond1}
        (\nabla^E\ot B)\circ \tau_{B,E} = (E\ot \widetilde{R})\circ (\tau_{B,E}\ot \O^1 A)\circ (B\ot \nabla^E),
    \end{equation}
which is the analogous of condition \eqref{eq:nabla2cond1}, and such that $(E\ot B)\oplus(A\ot F)$ becomes an $A\ot_R B$ bimodule with left action
    \[
        (\alpha\ot \beta)\cd (e\ot b,a\ot f):= (\alpha e_\tau\ot \beta_\tau b, \alpha a_R\ot \beta_R f),
    \]
then we have that
    \begin{eqnarray*}
        \nabla((\alpha\ot\beta)(e\ot b)) &=& \nabla_1(\alpha e_\tau \ot \beta_\tau b) = \\
        &=& (\alpha e_\tau)_i\ot 1\ot_{A\ot_R B} d_A(a'_i)\ot \beta_\tau b+\\
        && +\ \alpha e_\tau\ot 1\ot_{A\ot_R B} 1\ot d_B(\beta_\tau b) = \\
        &=& \alpha(e_\tau)_i\ot 1 \ot_{A\ot_R B} d_A(a_\tau)_i\ot\beta_\tau b+\\
        && +\ (e_\tau)_\vphi\ot 1\ot_{A\ot_R B} (d_A\alpha)_\vphi\ot \beta_\tau
b + \\
        && +\ \alpha e_\tau\ot 1\ot_{A\ot_R B} 1\ot d_B(\beta_\tau)b +\\
        && +\ \alpha
e_\tau\ot 1\ot_{A\ot_R B} 1\ot \beta_\tau d_B b
\overset{\eqref{eq:nabla1cond1}}{=} \\
        &\overset{\eqref{eq:nabla1cond1}}{=}&   \alpha(e_i)_\tau \ot 1
\ot_{A\ot_R B} (d_Aa_i)_{\widetilde{R}}\ot(\beta_\tau)_{\widetilde{R}} b+\\
        && \ \alpha e_\tau\ot 1\ot_{A\ot_R B} 1\ot \beta_\tau d_B b + \\
        &&+\ (e_\tau)_\vphi\ot 1\ot_{A\ot_R B} (d_A\alpha)_\vphi\ot \beta_\tau
b +\\
        && +\ \alpha e_\tau\ot 1\ot_{A\ot_R B} 1\ot d_B(\beta_\tau)b = \\
        &=& (\alpha\ot\beta)\nabla_1(e\ot b) +\\
        && +\  (e_\tau)_\vphi\ot 1\ot_{A\ot_R
B} (d_A\alpha)_\vphi\ot \beta_\tau b + \\
        &&+\ \alpha e_\tau\ot 1\ot_{A\ot_R B} 1\ot d_B(\beta_\tau)b.
    \end{eqnarray*}
On the other hand,
    \begin{eqnarray*}
        \nabla((\alpha\ot\beta)(a\ot f)) &=& \nabla_2(\alpha a_R \ot \beta_R f) = \\
        &=& 1\ot (\beta_R f)_\sigma \ot_{A\ot_R B} d_A((\alpha a_R)_\sigma)\ot 1
+\\
        && +\ \alpha a_R\ot (\beta_R f)_k\ot_{A\ot_R B}1\ot d_Bb'_k)
\overset{\eqref{eq:nabla2cond4}}{=} \\
        &\overset{\eqref{eq:nabla2cond4}}{=}& 1\ot (\beta_R f)_{\sigma\bar{\sigma}} \ot_{A\ot_R B} d_A (\alpha_{\bar{\sigma}}(a_R)_\sigma)\ot 1+ \\
        &&+\ \alpha a_R\ot (\beta_R f)_k\ot_{A\ot_R B} 1\ot d_B(b'_k)
\overset{\eqref{eq:nabla2cond5}}{=}\\
        &\overset{\eqref{eq:nabla2cond5}}{=}&  1\ot (\beta
f_\sigma)_{\bar{\sigma}} \ot_{A\ot_R B} d_A (\alpha_{\bar{\sigma}}a_\sigma)\ot 1
+\\
        && +\ \alpha a_R\ot (\beta_R f)_k\ot_{A\ot_R B} 1\ot d_B(b'_k) = \\
        &=& 1\ot (\beta f_\sigma)_{\bar{\sigma}} \ot_{A\ot_R B} d_A
(\alpha_{\bar{\sigma}})a_\sigma \ot 1 +\\
        && +\ 1\ot (\beta f_\sigma)_{\bar{\sigma}}
\ot_{A\ot_R B} \alpha_{\bar{\sigma}}d_A (a_\sigma)\ot 1 + \\
        && +\ \alpha a_R\ot \beta_R f_k\ot_{A\ot_R B} 1\ot d_Bb_k +\\
        && +\ \alpha a_R\ot f_\psi \ot_{A\ot_R B} 1\ot (d_B(\beta_R))_\psi =\\
        &=& (\alpha\ot \beta)\nabla_2(a\ot f) + \\
        && +\ 1\ot (\beta f_\sigma)_{\bar{\sigma}} \ot_{A\ot_R B} d_A
(\alpha_{\bar{\sigma}})a_\sigma \ot 1 + \\
        &&+\, \alpha a_R\ot f_\psi \ot_{A\ot_R B} 1\ot (d_B(\beta_R))_\psi.
    \end{eqnarray*}
Adding up these two equalities we  obtain
\begin{eqnarray*}
    \nabla((\alpha\ot\beta)(e\ot b, a\ot f)) &=& (\alpha\ot \beta)\nabla (e\ot
b, a\ot f) + \\
    &&+\ \xi(d(\alpha\ot\beta)\ot_{A\ot_R B} (e\ot b, a\ot f)),
\end{eqnarray*}
where the map $\xi: (\O^1A\ot B\oplus A\ot \O^1B)\ot_{A\ot_R B}
(E\ot B\oplus A\ot F) \to (E\ot B\oplus A\ot F)\ot_{A\ot_R
B}(\O^1A\ot B\oplus A\ot \O^1B)$ is defined by
$\xi:=\xi_{11}+\xi_{12}+\xi_{21}+\xi_{22}$, being
\begin{gather*}
    \xi_{11}(d_A\alpha\ot\beta\ot_{A\ot_R B} e\ot b) := (e_\tau)_\vphi\ot 1\ot_{A\ot_R B} (d_A\alpha)_\vphi\ot \beta_\tau
    b,\\
    \xi_{12}(\alpha\ot d_B\beta \ot_{A\ot_R B} e\ot b) := \alpha e_\tau\ot 1\ot_{A\ot_R B} 1\ot
    d_B(\beta_\tau)b,\\
    \xi_{21}(d_A\alpha\ot\beta\ot_{A\ot_R B} a\ot f) := 1\ot (\beta f_\sigma)_{\bar{\sigma}} \ot_{A\ot_R B} d_A (\alpha_{\bar{\sigma}})a_\sigma \ot
    1,\\
    \xi_{22}(\alpha\ot d_B \beta \ot_{A\ot_R B} a\ot f) := \alpha a_R\ot f_\psi \ot_{A\ot_R B} 1\ot (d_B(\beta_R))_\psi.
\end{gather*}
Hence, in order to show that the product connection $\nabla$ is a
bimodule connection we only have to show that $\xi$ is a bimodule
morphism, which is equivalent to prove that all the $\xi_{ij}$ are
bimodule morphisms.

\begin{lemma}\label{lemmaxi11}
    The map $\xi_{11}$ is a left $(A\ot_R B)$--module morphism, if, and only if, the equality
    \begin{equation}\label{xi11compatibility}
        (\vphi\ot B)\circ (\O^1 A\ot\tau_{B,E})\circ (\widetilde{R}\ot E)=(E\ot \widetilde{R})\circ (\tau_{B,E}\ot \O^1 A)\circ (B\ot \vphi)
    \end{equation}
    is satisfied in $B\ot \O^1 A\ot E$.
\end{lemma}

\begin{proof}
    In order to check that the compatibility condition is necessary, just apply the compatibility with the module action to an element of the form $1\ot b\ot \omega\ot 1\ot e\ot 1$.

    Conversely, assuming condition \eqref{xi11compatibility}, we have that
    \begin{multline*}
        \xi_{11}((x\ot y)\cd (d\alpha\ot\beta\ot_{A\ot_R B} e\ot b))
= \\
    = \xi_{11}(x(d\alpha)_{\widetilde{R}}\ot y_{\widetilde{R}}\beta \ot_{A\ot_R
B} e\ot b) = \\
        = (e_\tau)_\vphi\ot 1\ot_{A\ot_R B} (x(d\alpha)_{\widetilde{R}})_\vphi
\ot (y_{\widetilde{R}}\beta)_\tau b \overset{[1]}{=}\\
        \overset{[1]}{=} x(e_\tau)_\vphi\ot 1\ot_{A\ot_R B}
((d\alpha)_{\widetilde{R}})_\vphi \ot (y_{\widetilde{R}}\beta)_{\tau} b
\overset{[2]}{=}\\
        \overset{[2]}{=} x ((e_{\tau})_{\bar{\tau}})_\vphi\ot 1\ot_{A\ot_R B}
((d\alpha)_{\widetilde{R}})_{\vphi}\ot
(y_{\widetilde{R}})_{\bar{\tau}}\beta_\tau b
\overset{\eqref{xi11compatibility}}{=}\\
        \overset{\eqref{xi11compatibility}}{=}
x((e_\tau)_\vphi)_{\bar{\tau}}\ot 1\ot_{A\ot_R B}
((d\alpha)_{\vphi})_{\widetilde{R}}\ot
(y_{\bar{\tau}})_{\widetilde{R}}\beta_\tau b = \\
        = x((e_{\tau})_\vphi)_{\bar{\tau}}\ot 1\ot_{A\ot_R B} ((1\ot
y_{\bar{\tau}})\cd((d\alpha)_{\vphi}\ot \beta_\tau b)) = \\
        = x((e_{\tau})_\vphi)_{\bar{\tau}}\ot y_{\bar{\tau}}\ot_{A\ot_R B}
(d\alpha)_{\vphi}\ot \beta_\tau b =\\
        = (x\ot y)\xi_{11}(d\alpha\ot\beta\ot_{A\ot_R B} e\ot b),
    \end{multline*}
    where in [1] we are using that $\vphi$ is a left module map, in [2] that $\tau$ is a module twisting map.
\end{proof}
It is straightforward checking that $\xi_{11}$ is a right module
map, and thus left to the reader. In a completely analogous way, it
is straightforward to check that $\xi_{22}$ is a left module map,
whilst for the right module condition we need a compatibility
relation similar to \eqref{xi11compatibility}. More concretely, we
have the following result, whose  proof is analogous to the one of
Lemma \ref{lemmaxi11}:

\begin{lemma}
    The map $\xi_{22}$ is a right $A\ot_R B$--module morphism if, and only if, the equality
    \begin{equation}\label{xi22compatibility}
        (A\ot \psi)\circ(\widetilde{R}\ot F)\circ(\O^1 B\ot \tau_{F,A})=(\tau_{F,A}\ot \O^1 B)\circ(F\ot \widetilde{R})\circ (\psi\ot A)
    \end{equation}
    is satisfied in $\O^1 B\ot F\ot A$.
\end{lemma}

For $\xi_{12}$ and $\xi_{21}$, the right (resp. left) module map
conditions are also straightforward. We will show now that
$\xi_{12}$ is a left module map, the proof that $\xi_{21}$ is a
right module map being analogous.

    \begin{eqnarray*}
        \xi_{12}((x\ot y)\cd(\alpha\ot d\beta\ot_{A\ot_R B})e\ot b) &=& \xi_{12}(x\alpha_R\ot yd\beta \ot_{A\ot_R B} e\ot b) =\\
        &=& \xi_{12}(x\alpha_R\ot d(y_R\beta) \ot_{A\ot_R B} e\ot b) - \\
        && - \xi_{12}(x\alpha_R\ot d(y_R) \ot_{A\ot_R B} e_\tau \ot \beta_\tau b)=\\
        &=& x\alpha_R e_\tau \ot 1\ot_{A\ot_R B} 1\ot (d(y_R\beta)_\tau)b - \\
        && - x\alpha_R (e_\tau)_{\bar{\tau}}\ot 1\ot_{A\ot_R B} 1\ot d((y_R)_{\bar{\tau}}) \beta_{\tau} b\overset{[1]}{=}\\
        &\overset{[1]}{=}& x\alpha_R (e_\tau)_{\bar{\tau}} \ot 1\ot_{A\ot_R B} 1\ot (d((y_R)_{\bar{\tau}}\beta_\tau))b - \\
        && - x\alpha_R (e_\tau)_{\bar{\tau}}\ot 1\ot_{A\ot_R B} 1\ot d((y_R)_{\bar{\tau}}) \beta_{\tau} b =\\
        &=& x\alpha_R (e_\tau)_{\bar{\tau}} \ot 1\ot_{A\ot_R B} 1\ot d((y_R)_{\bar{\tau}})\beta_\tau b + \\
        && +x\alpha_R (e_\tau)_{\bar{\tau}} \ot 1\ot_{A\ot_R B} 1\ot (y_R)_{\bar{\tau}}(d(\beta_\tau))b-\\
        && - x\alpha_R (e_\tau)_{\bar{\tau}}\ot 1\ot_{A\ot_R B} 1\ot d((y_R)_{\bar{\tau}}) \beta_{\tau} b =\\
        &=& x\alpha_R (e_\tau)_{\bar{\tau}} \ot 1\ot_{A\ot_R B} 1\ot (y_R)_{\bar{\tau}}(d(\beta_\tau))b \overset{[2]}{=} \\
        &\overset{[2]}{=}& x(\alpha e_\tau)_{\bar{\tau}}\ot 1 \ot_{A\ot_R B} 1\ot y_{\bar{\tau}}(d\beta_\tau)b = \\
        &=& x(\alpha e_\tau)_{\bar{\tau}}\ot y_{\bar{\tau}} \ot_{A\ot_R B} 1\ot (d\beta_\tau)b = \\
        &=& (x\ot y)\cd \xi_{12}(\alpha\ot d\beta\ot_{A\ot_R B})e\ot b),
    \end{eqnarray*}
    where in [1] and [2] we use that $\tau_{F,A}$ is a module twisting map.

Summarizing, we have proved the following result:
\begin{theorem}
    Let $E$ be a bimodule over $A$, $(\nabla^E,\vphi)$ a bimodule connection on $E$, $F$ a bimodule over $B$, $(\nabla^F,\psi)$, $R:B\ot A\to A\ot B$ an invertible twisting map; $\tau_{F,A}:F\ot A\to A\ot F$ a right module twisting map satisfying condition \eqref{eq:nabla2cond1} and $\tau_{B,E}:B\ot E\to E\ot B$ a left module twisting map  satisfying condition \eqref{eq:nabla1cond1}. Assume also that conditions \eqref{xi11compatibility} and \eqref{xi22compatibility} are satisfied, then the product connection of $\nabla^E$ and $\nabla^F$ is a bimodule connection with respect to the morphism $\xi$.
\end{theorem}

\section{Examples} \label{examples}
\setcounter{equation}{0}

Let us start by recalling some facts from \cite{Cuntz95a}. For any
projective (right) module $E$ over an algebra $A$, there exists a
module $E'$ such that $E\oplus E' = A^n$, and we have two canonical
mappings
    \[
        p:A^n=E\oplus E' \lto E\quad \text{and } \lambda: E
        \hookrightarrow E\oplus E',
    \]
we can then define the map $\nabla_0:=(p\ot\Id)\circ (A^n\ot d)\circ
(\lambda\ot \Id)$ as the composition given by
    \[ \xymatrix{
            E\ot_A\O^p A \ar[rr]^{\lambda\ot \Id}&&
            A^n\ot_A \O^p A \ar[rr]^{A^n\ot d} &&
            \O^{p+1} A \ar[rr]^{p\ot \Id} && E\ot_A\O^{p+1}A
        }
    \]
The operator $\nabla_0$ is a (flat) connection on $E$, called the
\dtext{Grassmann connection} on $E$.

\begin{rem}
    Physicists sometimes use the shorthand notation $\nabla_0=pd$ to denote the
    Grassmann connection.
\end{rem}

It is also well known (cf. for instance \cite{Cuntz95a}) that the
space of all linear connections over a projective module $E$ is an
affine space modeled on the space of $A$--module morphisms
$\End_A(E)\ot_A \O^1 A$, and henceforth we can write any linear
connection $\nabla$ on $E$ as $\nabla=\nabla_0 + \alpha$, being
$\alpha\in \End_A(E)\ot_A \O^1 A$, where the ``\emph{matrix}''
$\alpha$ is called the \dtext{gauge potential} of the connection
$\nabla$.

\subsection{Product connections on the quantum plane $k_q[x,y]$}\hfill

Consider now $A:=k[x]$ the polynomial algebra in one variable. Since
any projective module over $A$ is free (actually, by Quillen-Suslin Theorem,
any projective module over any polynomial ring is free) it is enough to
consider connections for modules of the form $E=A^m$. If we denote by $\{e_i\}_{i=1,\dotsc,
m}$ the canonical generator set for $E$, we may write the Grassmann
connection on $E$ as
    \begin{equation}
        \nabla_0^E(a_1,\dotsc,a_m)=e_1\ot_A da_1 + \dotsb + e_m\ot_A
        da_m\in E\ot_A \O^1 A.
    \end{equation}
Analogously, let $B:=k[y]$, $F:=B^n$ with canonical generating
system $\{f_j\}_{j=1,\dotsc,n}$ and Grassmann connection
    \begin{equation}
        \nabla_0^F(b_1,\dotsc, b_n)=f_1\ot_B db_1 + \dotsb +
        f_n\ot_Bdb_n.
    \end{equation}

Recall that the quantum plane $k_q[x,y]$ may be seen as the twisted
tensor product $k[x]\ot_R k[y]$ with respect to the twisting map
obtained by extension of $R(y\ot x):=q x\ot y$. This is an
invertible twisting map which extends to an invertible module
twisting map $\tau_{F,A}:F\ot A\to A\ot F$ in a natural way. For elements $e\ot b \in E\ot B$, where $e=(a_1,\dotsc,a_m)$, and a generator $x\ot f$  with $f=(y^{i_1},\dotsc, y^{i_n})$ of $A\ot F$, using the definition of our product connection given by Equation \eqref{productconnection}, we have that the product of the Grassmann connections is
    \begin{eqnarray*}
        \nabla^{gr}(e\ot b,x\ot f) & = & \left( \sum e_i\ot 1\ot da_i\right)\ot b + e\ot 1\ot 1\ot db + \\
        & & + x\ot \left( \sum f_k\ot 1\ot dy^{i_k}\right) + 1\ot (q^{-i_1} y^{i_1},\dotsc,q^{-i_n} y^{i_n})\ot dx\ot 1
    \end{eqnarray*}

\begin{rem}
    If we introduce the notation $\lambda_q (p(y)):=p(qy)$, we can give the former expression for an element $a\ot f$ of the form $a=x^j$, $f=(b_1,\dotsc,b_n)\in F$ as
    \begin{eqnarray*}
        \nabla^{gr}(e\ot b,a\ot f) & = & \sum_i e_i\ot 1\ot da_i\ot b + e\ot 1\ot 1\ot db + \\
        & & + \sum_k a\ot f_k\ot 1\ot db_k + \sum_k 1\ot \lambda_{q^{-j}}(b_k)\ot d(x^j)\ot 1
    \end{eqnarray*}
\end{rem}

Now, for a generic connection $\nabla^E$ over the module $E$, there must exist
a potential $\alpha^E=\vphi_i\ot \omega_i\in
\End{E}\ot_A\O^1 A$ given by $\alpha^E(a_1,\dotsc,a_m)=\sum_{i,j}
\vphi_i(a_j)\ot \omega_i$ such that $\nabla^E=\nabla^E_0+\alpha^E$. In the same
way, for a generic connection $\nabla^F$ on $F$ there must exist a potential
$\alpha^F=\sum_k \psi_k\ot \eta_k$, given by $\alpha^F(b_1,\dotsc,
b_n)=\sum_{k,l}\psi_k{b_l}\ot\eta_k$, and such that
$\nabla^F=\nabla^F_0+\alpha^F$. Applying the formula for the product connection
to $\nabla^E$ and $\nabla^F$ we easily observe that
    \begin{equation*}
        \nabla(e\ot b,a\ot f) = \nabla^{gr}(e\ot b,a\ot f) + \sum_{i,j}\vphi_i(a_j)\ot 1\ot \omega_i\ot b + \sum_{k,l}a\ot\psi_k(b_l)\ot 1\ot\eta_k,
    \end{equation*}
expression that tells us the formula for all possible product connections on the quantum plane.

 \nocite{Landi97a} \nocite{Gracia01a} \nocite{Mourad95a}
 \nocite{BeggsUNa} \nocite{Jara03a} \nocite{DuboisViolette95a}
 \nocite{Connes86a} \nocite{LopezUNa} \nocite{Beck69a}
\providecommand{\bysame}{\leavevmode\hbox
to3em{\hrulefill}\thinspace}
\providecommand{\MR}{\relax\ifhmode\unskip\space\fi MR }
\providecommand{\MRhref}[2]{%
  \href{http://www.ams.org/mathscinet-getitem?mr=#1}{#2}
} \providecommand{\href}[2]{#2}

\end{document}